\documentclass[12pt,reqno]{elsarticle}

\usepackage[usenames]{color}
\usepackage{amssymb}
\usepackage{graphicx}
\usepackage{amscd}

\usepackage[colorlinks=true,
linkcolor=webgreen,
filecolor=webbrown,
citecolor=webgreen]{hyperref}

\definecolor{webgreen}{rgb}{0,.5,0}
\definecolor{webbrown}{rgb}{.6,0,0}

\usepackage{color}
\usepackage{fullpage}
\usepackage{float}

\usepackage{graphics,amsmath,amssymb}
\usepackage{amsthm}
\usepackage{amsfonts}
\usepackage{latexsym}
\usepackage{epsf}
\usepackage{subfigure}

\begin{document}

\begin{center}
\vskip 1cm{\LARGE\bf Counting symmetry classes of dissections of a convex regular polygon
}
\vskip 1cm
\large
Douglas Bowman and Alon Regev\\
Department of Mathematical Sciences\\
Northern Illinois Univeristy\\
DeKalb, IL
\end{center}

{\center \section*{Abstract}}
This paper proves explicit formulas for the number of dissections of a convex regular polygon modulo the action of the cyclic and dihedral groups. The formulas are obtained by making use of the Cauchy-Frobenius Lemma as well as bijections between rotationally symmetric dissections and simpler classes of dissections. A number of special cases of these formulas are studied. Consequently, some known enumerations are recovered and several new ones are provided.

\vskip .2 in

\theoremstyle{plain}
\newtheorem{theorem}{Theorem}
\newtheorem{corollary}[theorem]{Corollary}
\newtheorem{lemma}[theorem]{Lemma}
\newtheorem{proposition}[theorem]{Proposition}

\theoremstyle{definition}
\newtheorem{definition}[theorem]{Definition}

\newtheorem{remark}[theorem]{Remark}

\newcommand{\PP}{\mathcal{C}}

\newcommand{\summ}{\sum\limits}
\newcommand{\prodd}{\prod\limits}

\newcommand{\G}[2]{A^{#1}_{#2}}

\newcommand{\R}{\rho}
\newcommand{\T}{\tau}
\newcommand{\E}{\varepsilon}

\newcommand {\hh}{\gamma}

\section{Introduction}

\begin{figure}\label{P3}
\begin{center}
\epsfxsize=3.5in
 \includegraphics[scale =0.8] {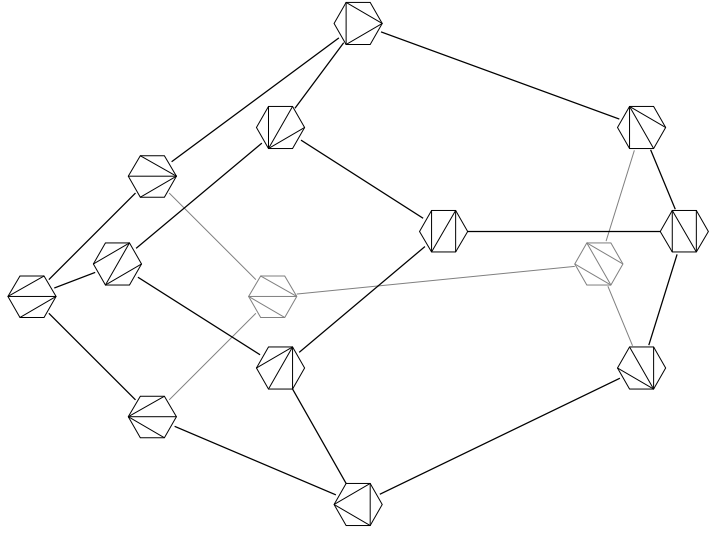}
 \end{center}
\caption{The three-dimensional associahedron}
\end{figure}

In 1963 Moon and Moser \cite{MM} enumerated the equivalence classes of triangulations of a regular convex $n$-gon modulo the action of the dihedral group $D_{2n}$. A year later, Brown \cite{Br} enumerated the equivalence classes of these triangulations modulo the action of the cyclic group $Z_n$.
Recall that the triangulations of an $n$-gon are in bijection with the vertices of the associahedron of dimension $n-3$ (see Figure \ref{P3}).
Lee \cite{Le} showed that the associahedron can be realized as a polytope in $(n-3)$-dimensional space having the dihedral symmetry group $D_{2n}$. Thus Moon and Moser's result and Brown's result are equivalent to enumerating of the vertices of the associahedron modulo the dihedral action and the cyclic action, respectively. The enumeration by
Moon and Moser also arose recently in the work of  Ceballos, Santos and Ziegler \cite{CSZ}. Their work  describes a family of realizations  of the associahedron (due to Santos), and proves that the number of normally non-isomorphic realizations is the number of triangulations of a regular polygon modulo the dihedral action. In this paper we generalize the
results of Moon and Moser, as well as Brown, and enumerate the number of {\em dissections} of regular polygons modulo the dihedral and cyclic actions.

\begin{definition}
Let $n\ge 3$. A {\em $k$-dissection} of an $n$-gon is a partition of the $n$-gon into $k+1$ polygons by $k$ non-crossing diagonals. A {\em triangulation} is an $(n-3)$-dissection of an $n$-gon and an {\em almost-triangulation} is an $(n-4)$-dissection. Let $G(n,k)$ be the set of $k$-dissections of an $n$-gon, and let $G(n)=\bigcup\limits_{k=0}^{n-3}G(n,k)$.
\end{definition}

In terms of associahedra, a $k$-dissection corresponds to an $(n-k-3)$-dimensional face on an associahedron of dimension $n$. A natural generalization of the results of Moon and Moser and of Brown is the enumeration of $G(n,k)/D_{2n}$ and $G(n,k)/Z_n$, the sets of cyclic and dihedral classes, respectively, in $G(n,k)$. In 1978, Read \cite{R} considered an equivalent problem. He enumerated certain classes of cellular structures, which are in bijection with $G(n,k)/D_{2n}$ and
$G(n,k)/Z_n$. Read found generating functions for the number of such classes, and included tables of values \cite[Tables 3 and 5]{R}. In fact, the first diagonal of Table 5 of Read corresponds to the sequence found by Moon and Moser, and the first diagonal of Table 3 of Read corresponds to the sequence found by Brown. Lisonek \cite{Li} studied these results of Read 
and showed that the sequences $|G(n,k)/D_{2n}|$ and $|G(n,k)/Z_n|$ are ``quasi-polynomial" in $n$ when $k$ is fixed. (Here and throughout, $|X|$ denotes the cardinality of a finite set $X$.) More recently, Read and Devadoss \cite{DR} studied various equivalence relations on the set of polygonal dissections. They gave a sequence of figures
\cite[Figures 22-25]{DR} representing all the dihedral classes of $n$-gons for $3 \le n\le 9$. However, none of the above authors give an explicit formula for $|G(n,k)/D_{2n}|$ or for $|G(n,k)/Z_n|$.
The present authors \cite{BR} give an explicit formula enumerating $G(n,n-4)/D_{2n}$, the dihedral classes of almost-triangulations, equivalently, of edges of associahedra. This formula agrees with the values of the second diagonal of Table 5 of Read \cite{R}.

Explicit formulas for $|G(n,k)/Z_n|$ and $|G(n,k)/D_{2n}|$ could in principle be derived from Read's iteratively defined generating functions, but the resulting formulas would be considerably more complicated than those computed here; see equations \eqref{diheq} and \eqref{cyceq}.
Our approach to solving these enumeration problems is similar to that of Moon and Moser \cite{MM}. For each element of the dihedral group, the number of dissections in $G(n,k)$ which are fixed under its action is computed. The Cauchy-Frobenius Lemma is then used to derive the number of dihedral and cyclic classes in $G(n,k)$.

In Section \ref{rotsec} we introduce a combinatorial bijection \eqref{bijeq} between certain rotationally symmetric dissections ({\em centrally unbordered} dissections, see Definition \ref{classes}) and a set  $G^*(n,k)$ of {\em marked dissections}, which are dissections with one of their parts distinguished. These marked dissections are easy to generate and enumerate. A bijection for {\em centrally bordered} dissections is implicit in the proof of Lemma \ref{bdd}.
 Przytycki and Sikora \cite{PS} studied a set of marked dissections $P_i(s,n)$, which is a subset of $G^*(n,k)$; however, the classes of dissections enumerated in \cite{PS}  are different from those studied here.

Besides their intrinsic interest, bijections involving polygonal dissections have connections to other mathematical structures; for example, Torkildsen \cite{T} proved a bijection between $G(n,n-3)/Z_{n}$ and the mutation class of quivers of Dynkin type $A_n$, while Przytycki and Sikora describe a relationship between their bijection (between $P_i(s,n)$ and another combinatorial structure) and their work in knot theory as well as Jones' work on planar algebras.

After proving the general formulas \eqref{diheq} and \eqref{cyceq} in Sections  \ref{Prelim} through \ref{rotsec},  special cases are studied in Section \ref{specs}. Consequently we not only recover known enumerations but are able to provide several that are new.  We note several of the interesting special cases here. For example, setting $k=n-3$ in \eqref{diheq} recovers the result of Moon and Moser \cite{MM}, and setting $k=n-3$ in \eqref{cyceq} recovers the result of Brown \cite{Br}.
Setting $k=n-4$ in \eqref{diheq} recovers the result of the authors \cite{BR}, while the following theorem gives a formula for the number of cyclic classes in the case $k=n-4$.
For a nonnegative integer $n$, let $C_n$ denote the $n$-th Catalan number; \[C_{n}={1\over n+1}{2n\choose n},\] and let $C_n=0$ otherwise.

\begin{theorem}\label{cycn4}
Let $n\ge 4$. The number $|G(n,n-4)/Z_n|$ of almost-triangulations of an $n$-gon (equivalently, edges of the $(n-3)$-dimensional associahedron) modulo the cyclic action is given by
\begin{equation}\label{cycn4eq}
\frac{n-3}{2n}C_{n-2}+\frac{1}{2}C_{n/4-1}+\frac{1}{4}C_{n/2-1}.
\end{equation}
\end{theorem}

Setting $k=n-5$ gives the following formulas.

\begin{theorem}\label{k=n-5}
Let $n\ge 5$.
\begin{enumerate}
\item The number $|G(n,n-5)/Z_n|$ of $(n-5)$-dissections of an $n$-gon (equivalently, the number of two-dimensional faces of the $(n-3)$-dimensional associahedron) modulo the cyclic action is given by
\begin{align}\label{cycn5eq}
{(n-3)^2(n-4)\over 4n(2n-5)}C_{n-2}+{n-4\over 8}C_{n/2-1} + {4 \over 5}C_{n/5-1}.
\end{align}
\item The number $|G(n,n-5)/D_{2n}|$ of $(n-5)$-dissections of an $n$-gon (equivalently, the number of two-dimensional faces of the $(n-3)$-dimensional associahedron) modulo the dihedral action is given by
\begin{align*}
{(n-3)^2(n-4)\over 8n(2n-5)}C_{n-2}+{2\over 5}C_{n/5-1}+{3(n-4)(n-1)\over 16(n-3)}C_{n/2-1},
\end{align*}
if $n$ is even, and
\begin{align}\label{dihn5eq}
{(n-3)^2(n-4)\over 8n(2n-5)}C_{n-2}+{2\over 5}C_{n/5-1}+{n^2-2n-11 \over 8(n-4)}C_{(n-3)/2},
\end{align}
if $n$ is odd.
\end{enumerate}
\end{theorem}
A formula equivalent to \eqref{cycn4eq} occurs in the Online Encyclopedia of Integer Sequences \cite[sequence A0003444]{S}, while the sequences of \eqref{cycn5eq} and \eqref{dihn5eq} occur there without a formula \cite[sequences A0003450 and A0003445]{S}.

Finally, setting $k=n-6$ in \eqref{cyceq} gives the following formula.

\begin{theorem}\label{k=n-6}
Let $n\ge 6$.
The number $|G(n,n-6)/Z_n|$ of $(n-6)$-dissections of an $n$-gon (equivalently, the number of three-dimensional faces of the $(n-3)$-dimensional associahedron) modulo the cyclic action is given by
\[{(n-3)(n-4)^2(n-5)\over 24n(2n-5)}C_{n-2}+{(n-4)^2\over 4n}C_{n/2-2}+{n-3\over 9}C_{n/3-1}+{1\over 3} C_{n/6-1}.\]
\end{theorem}

Another set of results is obtained by specializing equations \eqref{diheq} and \eqref{cyceq} to fixed values of $k$. Setting $k=1$ gives the formulas $|G(n,1)/Z_n|=|G(n,1)/D_{2n}|=n/2-1$ if $n$ is even and
$|G(n,k)/Z_n|=|G(n,1)/D_{2n}|=(n-3)/2$ if $n$ is odd: these formulas are easy to see directly. In the case $k=2$, \eqref{diheq} and \eqref{cyceq} are more interesting.

\begin{theorem}\label{k=2}
Let $n\ge 2$.
\begin{enumerate}
\item The number of $2$-dissections of an $n$-gon modulo the cyclic action is
\[|G(n,2)/Z_n|=\begin{cases} {1\over 12}n(n-2)(n-4), \qquad \text{ if $n$ is even}\\
{1\over 12}(n+1)(n-3)(n-4), \qquad \text{ if $n$ is odd.}\end{cases}\]
\item The number of $2$-dissections of an $n$-gon modulo the dihedral action is
\[|G(n,2)/D_{2n}|=\begin{cases}
{1\over 24}(n-4)(n-2)(n+3),  \qquad \text{ if $n$ is even}\\
{1\over 24}(n-3)(n^2-13), \qquad \text{ if $n$ is odd.}
\end{cases}\]
\end{enumerate}
\end{theorem}
Note that Theorem \ref{k=2} agrees with the result of Lisonek \cite{Li} in that the formulas obtained are quasi-polynomials.

\subsection{The general formulas}

Let $\G{n}{k}=|G(n,k)|$ be the number of $k$-dissections of an $n$-gon.
Cayley \cite{C} showed that for integers $0 \le k\le n-3$, 
\begin{equation}\label{disseq}
\G{n}{k}=
\frac{1}{k+1}{n+k-1 \choose k}{n-3 \choose k}.
\end{equation}
We take $\G{2}{0}=1$ corresponding to the trivial dissection of a digon ($2$-gon). Otherwise, unless $n$ and $k$ are integers with $0\le k \le n-3$,
let $\G{n}{k}=0$.
Note that
\begin{equation}\label{catdiss}
\G{n}{n-3}=\begin{cases} 0, \text{ if } n=2\\ C_{n-2}, \text{ otherwise.}\\ \end{cases}
\end{equation}

Let $\varphi(n)$ denote Euler's totient function (the number of positive integers less than $n$ that are relatively prime to $n$). The following two theorems are the main results of this paper.

\begin{theorem}\label{dih}
Let $1\le k \le n-3$. Let $|G(n,k)/D_{2n}|$ be the number of $k$-dissections of an $n$-gon (equivalently, $(n-k-3)$-dimensional faces of the $(n-3)$-dimensional associahedron) modulo the dihedral action.
If $n$ is even then $|G(n,k)/D_{2n}|$ is given by
\begin{align*}
    &{1\over 2n}\G{n}{k}        
    +{1\over 2} \G{n/2+1}{(k-1)/2}  
    + {1\over 4} \G{n/2+1}{k/2} 
    \nonumber  \\ &  
    + \summ_{3\le d|n}{\varphi(d)\over 2d} \G{n/d+1}{k/d-1} 
    + \summ_{2\le d\le n/3}{\varphi(d)(n+k-d) \over 2dn}\G{n/d}{k/d-1}  \nonumber \\         
    &+ \summ_{\substack{2\le d|n; \ r\ge 3; \ n_1+\ldots+n_r=n/d\\ k_1+ \ldots+k_r + |\{i : n_i\ge 2\}|= k/d}} {\varphi(d)\over 2r}  \prodd_{i=1}^r\G{n_i+1}{k_i} \nonumber \\ &
         + {1\over 4}  \summ_{\substack{1\le t \le k\\ n_0+\ldots + n_t=n/2\\k_0+\ldots + k_t=(k-t)/2}}  \ \
             \prodd_{s=0}^ {t} \big (\G{n_s+1}{k_s-1} + \G{n_s+1}{k_s} \big) \nonumber \\    
   & + {1\over 4} \summ_{\substack{0\le t \le k\\ n_0+\ldots + n_t=n/2-1\\ k_0+\ldots + k_t=(k-t)/2}} \ \
        \prodd_{s=0}^ {t} \big (\G{n_s+1}{k_s-1} + \G{n_s+1}{k_s} \big),                    
\end{align*}
and if $n$ is odd then $|G(n,k)/D_{2n}|$ is given by
\begin{align}\label{diheq}
   & {1\over 2n}\G{n}{k}        
    + \summ_{3\le d|n}{\varphi(d)\over 2d} \G{n/d+1}{k/d-1} 
    + \summ_{2\le d\le n/3}{\varphi(d)(n+k-d) \over 2dn}\G{n/d}{k/d-1}  \nonumber \\         
    &+ \summ_{\substack{2\le d|n; \ r\ge 3; \ n_1+\ldots+n_r=n/d\\ k_1+ \ldots+k_r + |\{i : n_i\ge 2\}|= k/d}} {\varphi(d)\over 2r}  \prodd_{i=1}^r\G{n_i+1}{k_i} \nonumber \\ &
         + {1\over 2} \summ_{\substack{1\le t \le k\\ n_0+\ldots + n_t=n/2\\k_0+\ldots + k_t=(k-t)/2}}  \ \
             \prodd_{s=0}^ {t} \big (\G{n_s+1}{k_s-1} + \G{n_s+1}{k_s} \big). \nonumber \\    
\end{align}
\end{theorem}

\begin{theorem}\label{cyc}
Let $n\ge 3$. The number $|G(n,k)/Z_n| $ of $k$-dissections of an $n$-gon (equivalently, ($n-k-3$)-dimensional faces of the $(n-3)$-dimensional associahedron) modulo the cyclic action is given by
\begin{align}\label{cyceq}
    {1\over n}\G{n}{k}        
    + \summ_{3\le d|n}{\varphi(d)\over d} \G{n/d+1}{k/d-1} 
    + \summ_{2\le d\le n/3}{\varphi(d)(n+k-d) \over dn}\G{n/d}{k/d-1}  \nonumber \\         
    + \summ_{\substack{2\le d|n; \ r\ge 3; \ n_1+\ldots+n_r=n/d\\ k_1+ \ldots+k_r + |\{i : n_i\ge 2\}|= k/d}} {\varphi(d)\over r}  \prodd_{i=1}^r\G{n_i+1}{k_i}. \nonumber \\ &
\end{align}
\end{theorem}

\section{Preliminaries}\label{Prelim}

For any labeled graph $H$, let $V(H)$ be its set of vertices and let $E(H)$ be its set of edges, defined to be two-element subsets of $V(H)$.
We frequently denote the edge $\{x,y\}$ by $xy$.
To a dissection $\Phi\in G(n)$ we associate a labeled graph $(V(\Phi),E(\Phi))$, where $V(\Phi)=\{0,\ldots, n-1\}$. To the sides of the $n$-gon we associate
the edges $S(\Phi)=\{01, 12, 23,\ldots ,(n-2)(n-1), (n-1)0\}$ and to the diagonals of the dissection we associate the rest of the edges of the graph. The {\em distance} between any two vertices $x,y\in V(\Phi)$ is defined to be the graph-theoretic distance between them in the subgraph $(V(\Phi),S(\Phi))$.

It is easily seen that two crossing diagonals of a convex $n$-gon correspond to edges $ab$ and $cd$ (with $a<b$ and $c<d$) if and only if
\begin{equation}\label{crossing}
a< c < b <d \, \, \text{ or } \, \, c < a < d < b.
\end{equation}
Thus if $\Phi\in G(n)$, there are no edges $ab, cd\in E(\Phi)$ satisfying \eqref{crossing}.
Hereafter we identify dissections with their labeled graphs.

The elements of the dihedral group $D_{2n}$ are denoted using $\E$, $\R$ and $\T$ to represent the identity, rotation by $2\pi/n$ and
reflection about a symmetry axis (which by convention passes through one of the vertices of the $n$-gon), respectively.
The cyclic group $Z_n$ can be identified with the subgroup of $D_{2n}$ generated by $\R$.

We use the notation $[x]_n$ to denote the remainder when $x$ is divided by $n$. The elements of $D_{2n}$ can be represented by their action on the vertices,
$\R(v)=[v+1]_n$ and $\T(v)=[-v]_n$.

For any $\sigma\in D_{2n}$, let $G(n,k; \sigma)$ denote the subset of $G(n,k)$ consisting of dissections fixed under the action of $\sigma$, and let
\[G(n;\sigma)=\bigcup\limits_{ k= 0}^{n-3}G(n,k;\sigma).\] Thus if $\Phi\in G(n)$ then $\Phi\in G(n;\R^i)$  if and only if
\begin{equation*}
 [x+i]_n[y+i]_n\in E(\Phi) \text{ whenever } [x]_n[y]_n\in E(\Phi),
\end{equation*}
and $\Phi\in G(n;\T\R^i)$  if and only if
\begin{equation*}
[i-x]_n[i-y]_n\in E(\Phi) \text{ whenever } [x]_n[y]_n\in E(\Phi) .
\end{equation*}

The Cauchy-Frobenius Lemma \cite{B} gives the equations
\begin{equation}\label{dihgeneq}
|G(n,k)/D_{2n}| = {1 \over 2n}\left( \summ_{i=0}^{n-1} |G(n,k; \R^i)| + \summ_{i=0}^{n-1} |G(n,k; \T \R^i)|\right)
\end{equation}
and
\begin{equation}\label{cycgeneq}
|G(n,k)/Z_n| = {1 \over n} \summ_{i=0}^{n-1} |G(n,k; \R^i)|.
\end{equation}
Equations \eqref{dihgeneq} and \eqref{cycgeneq} reduce the problem of enumerating the dihedral and cyclic classes in $G(n,k)$ to that of finding $G(n,k; \sigma)$ for each $\sigma\in D_{2n}$.
In fact, the following lemma shows that it suffices to consider only a subset of $D_{2n}$.
\begin{lemma}
Let $0\le i \le n-1$. Then
\begin{enumerate}
\item
\begin{equation*}
|G(n,k;\R^i)|=|G(n,k;\R^{\gcd(n,i)})|.
\end{equation*}
\item 
\begin{enumerate}
\item If $n$ and $i$ are even then $|G(n,k;\T\R^i)|=|G(n,k;\T)|$.
\item If $n$ is even and $i$ odd then $|G(n,k;\T \R^i)|=|G(n,k;\T \R)|$.
\item If $n$ is odd then $|G(n,k;\T \R^i)|=|G(n,k;\T)|$.
\end{enumerate}
\end{enumerate}
\begin{proof}
\begin{enumerate}
\item Since $\R^i$ and $\R^{\gcd(n,i)}$ generate the same subgroup of $D_{2n}$, they fix precisely the same elements of $G(n,k)$.
\item It is well known \cite[p. 243]{A} that conjugate elements in a group acting on a set have the same number of fixed points.
The results then follow
from the conjugacy relations in $D_{2n}$.
\end{enumerate}
\end{proof}
\end{lemma}

Thus equations \eqref{dihgeneq} and \eqref{cycgeneq} imply
\begin{equation}\label{dihgeneq2}
|G(n,k)/D_{2n}| = \begin{cases} {1 \over 2n}\summ_{d|n} \varphi(d)|G(n,k; \R^{n/d})|+{1\over 4} |G(n,k; \T)|+{1\over 4} |G(n,k; \T \R)|, & \text{ if $n$ is even} \\
{1 \over 2n}\summ_{d|n} \varphi(d)|G(n,k; \R^{n/d})|+{1\over 2} |G(n,k; \T)|, & \text{ if $n$ is odd}
\end{cases}
\end{equation}
and
\begin{equation}\label{cycgeneq2}
|G(n,k)/Z_n| = {1 \over n} \summ_{d|n} \varphi(d)|G(n,k; \R^{n/d})|.
\end{equation}
Theorems \ref{dih} and \ref{cyc} follow from calculating the terms in \eqref{dihgeneq2} and \eqref{cycgeneq2}, respectively.
Section \ref{axsec} addresses the terms $|G(n,k; \T)|$ and $|G(n,k; \T \R)|$, and Section \ref{rotsec} addresses the terms $|G(n,k;\R^{n/d})|$.

\section{Axially symmetric dissections}\label{axsec}

The sets $G(n,k;\T)$ and $G(n,k;\T \R)$ of axially symmetric dissections can be enumerated by considering the number of {\em perpendiculars}, i.e., diagonals of a dissection which are perpendicular to the axis of symmetry; these diagonals have the form $[v]_n[-v]_n$. Denote by $G(n,k;\T; t)$ the set of dissections in $G(n,k;\T)$ with exactly $t$ perpendiculars.
The notation $G(n,k;\T \R;t)$ is defined analogously. Thus
\[G(n,k;\T)=\sum_{t\ge 0}G(n,k;\T;t),\]
with the analogous formula holding for $G(n,k;\T \R;t)$.

\begin{figure}
\begin{center}
\epsfxsize=2.5in
\epsfbox{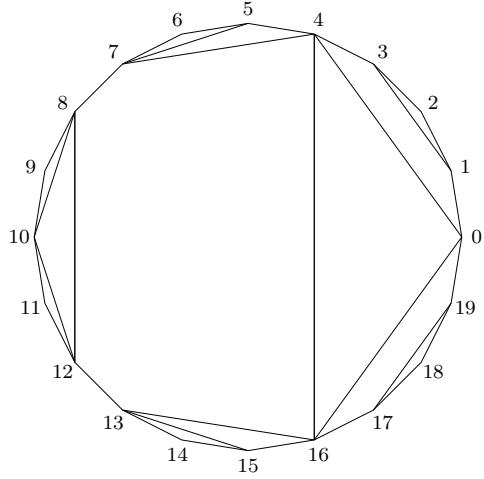}
\end{center}
\caption{An axially symmetric dissection.}
\label{20gon}
\end{figure}

\begin{lemma}\label{Tlem}
If $n$ is even then
\begin{equation}\label{p=0}
|G(n,k;\T;0)|= \G{n/2+1}{(k-1)/2}+\G{n/2+1}{k/2}
\end{equation}
and for $t\ge 1$,
\begin{equation}\label{p>0}
|G(n,k;\T;t)|= \summ_{\substack{n_0+\ldots + n_t=n/2\\k_0+\ldots + k_t=(k-t)/2}}
 \prodd_{s=0}^ {t} (\G{n_s+1}{k_s-1} + \G{n_s+1}{ k_s}).
\end{equation}

\begin{proof}
Let $[v_1]_n[-v_1]_n,\ldots ,[v_t]_n[-v_t]_n\in E(\Phi)$
be the perpendiculars of $\Phi$, where $t\ge 0$ and $0=v_0<v_1<\ldots <v_t<v_{t+1}=n/2$. Let  $n_s=v_{s+1}-v_s$ for $s=0,\ldots, t$. The case $t=0$ is considered separately since in this case $v_0v_1=0{n \over 2}=[-0]_n[-{n\over 2}]_n=[-v_0]_n[-v_1]_n$, while in all other cases
$v_sv_{s+1}\ne [-v_s]_n[-v_{s+1}]_n$.
Let $\Phi\in G(n,k;\T;0)$ and suppose first that $0{n\over 2}\in E(\Phi)$.
The remaining $k-1$ diagonals of $\Phi$ are equally distributed between the two sides of the symmetry axis. Each such dissection then
uniquely corresponds to a dissection of the resulting $(n/2+1)$-gon on either of its sides. Thus there are $\G{n/2+1}{(k-1)/2}$ dissections of this type.
By a similar argument, there are $\G{n/2+1}{k/2}$ dissections in $G(n,k;\T;0)$ which do not contain the diagonal $0{n\over 2}$, and \eqref{p=0} follows.
Now suppose $1\le t\le k$ and let $\Phi\in G(n,k;\T;t)$.
The $k-t$ other diagonals of $\Phi$ are pairs of the form $xy$ and $[-x]_n[-y]_n$, where $v_s\le x<y\le v_{s+1}$ and $0\le s \le t$.
Each such diagonal $xy$ is either of the form $v_sv_{s+1}$ or it is a diagonal of the $(n_s+1)$-gon with vertices $v_s, v_s+1, \ldots, v_{s+1}$.
Let $k_s$ be half the number of diagonals in the region defined by the vertices $v_s$, $v_{s+1}$, $[-v_{s+1}]_n$ and $[-v_s]_n$. If $v_sv_{s+1}\in E(\Phi)$ then the dissection of this region
corresponds to a $(k_s-1)$-dissection of the $(n_s+1)$-gon. Otherwise, it corresponds to a $k_s$-dissection
of the $(n_s+1)$-gon. This proves \eqref{p>0}.
\end{proof}
\end{lemma}
Figure \ref{20gon} gives an example of an axially symmetric dissection where, using the notation above, $n=20$, $k=12$, $t=2$, $v_1=4$, $v_2=8$, $k_0=2$, $k_1=2$ and $k_2=1$. The dissection of the region $s=0$ corresponds to a $1$-dissection of the pentagon with vertices $0, 1, 2, 3, 4$, and the dissection of the region $s=1$ corresponds to a $2$-dissection of the pentagon with vertices $4,5,6,7,8$.

For the next two lemmas,  $v_sv_{s+1}\ne [-v_s]_n[-v_{s+1}]_n$ for all $s$, and therefore the case $t=0$ need not be considered separately. The proofs are otherwise analogous to that of Lemma \ref{Tlem}.

\begin{lemma}\label{Tlem-odd}
If $n$ is even then for any $t\ge 0$
\begin{equation}\label{TReq}
 |G(n,k;\T \R;t)|  =
\summ_{\substack{n_0+\ldots + n_t=n/2-1\\k_0+\ldots + k_t=(k-t)/2}}
\ \prodd_{s=0}^ {t} (\G{n_s+1}{k_s-1} + \G{n_s+1}{k_s}).
\end{equation}
\end{lemma}

\begin{lemma}\label{TRlem}
If $n$ is odd then for $t\ge 0$,
\begin{equation}\label{Todd}
|G(n,k;\T;t)| =
\summ_{\substack{n_0+\ldots + n_t=(n-1)/2\\k_0+\ldots + k_t=(k-t)/2}}
\ \prodd_{s=0}^ {t} (\G{n_s+1}{k_s-1}+\G{n_s+1}{k_s}).
\end{equation}
\end{lemma}

\section{Components and marked dissections}

A dissection $\Phi\in G(n)$ can be associated with the set $\PP(\Phi)$ of {\em components} comprising it, each of these components being a polygon free of dissecting diagonals.
Thus a component is a subgraph $\gamma$ of $\Phi$ such that for some $r\ge 2$ and $0\le v_r=v_0<\ldots <v_{r-1}\le n-1$,

\[V(\gamma)=\{v_0,\ldots, v_{r-1}\},\]
\[E(\gamma)=\{v_0v_1,\ldots  , v_{r-1}v_r\},\]
and
\begin{align}\label{component}
v_iv_j\in E(\Phi) \implies v_iv_j\in E(\gamma).
\end{align}

 For example, if $\Phi$ is the dissection shown in Figure \ref{20gon} then $\PP(\Phi)$ consists of $32$ digons, $10$ triangles, two quadrilaterals and one hexagon. For $r\ge 2$ let $\PP_r(\Phi)$ be the set of $r$-gons in $\PP(\Phi)$, which we call {\em $r$-components}.

In what follows, $X$ represents any list of parameters. For any $r\ge 2$, an {\em $r$-marked dissection} is a dissection $\Phi$ with one of its $r$-components $\hh$ distinguished.
Let \[G^r(X)=\{(\Phi, \hh): \Phi \in G(X), \hh \in \PP_r(\Phi) \}\] be the set of  $r$-marked dissections associated with $G(X)$, and
let $G^*(X)=\bigcup\limits_{r\ge 2}G^r(X)$.

\begin{lemma}\label{marked}
Let $0\le k \le n- 3$. Then
\begin{equation}\label{r=2}
|G^2(n,k)|=(n+k)\G{n}{k},
\end{equation}
 and for $r \ge 3,$
\begin{equation}\label{MD}
r|G^r(n,k)|=n \summ_{\substack{n_1+\ldots+n_r=n\\ k_1+\ldots +k_r +|\{i:n_i\ge 2\}|= k}} \ \ \prodd_{i=1}^r\G{n_i+1}{k_i}.
\end{equation}
\begin{proof}
Equation \eqref{r=2} holds since $|\PP_2(\Phi)|=n+k$  for any dissection $\Phi\in G(n,k)$.
Let $r\ge 3$.
The left-hand side of \eqref{MD} enumerates the marked dissections having one of the $r$ vertices $v_0$ in the distinguished $r$-gon distinguished.
The right-hand side enumerates the same elements by selecting the $r$-component first and then dissecting the region between each side of the component and the $n$-gon. Choose $0\le v_0\le n-1$. Decompose the cycle $v_0[v_0+1]_n, [v_0+1]_n[v_0+2]_n,\ldots , [v_0+n-1]_n[v_0+n]_n$
into consecutive paths of length $n_i$ (where $1\le i\le r$ and $n_i\ge 1$) with
vertices $v_0,v_1,\ldots , v_r$.
Observe that every such decomposition of the edges of the $n$-gon corresponds to a set of $n_i$ satisfying $n_1+\ldots +n_r=n$. The vertices $V(\hh)=\{v_0, \ldots,v_{r-1}\}$ of an $r$-component are thus determined by
\[v_i=[v_0+\summ_{j=1}^i n_j]_n.\]
Now for each  $0\le i\le r-1$ select a dissection of the region between the edge $v_iv_{i+1}$ and the path from $v_i$ to $v_{i+1}$ along the sides of the $n$-gon. Each such dissection corresponds to a set of $k_i$ satisfying $k_1+\ldots +k_r +|\{i:n_i\ge 2\}|= k$. (The term $|\{i:n_i\ge 2\}|$ accounts for those sides of the $r$-component that are not sides of the $n$-gon). Finally for each $k_i$ there $\G{n_i+1}{k_i}$ such dissections.
\end{proof}
\end{lemma}

Since a triangulation consists of $n-2$ triangles, a simpler formula for $3$-marked triangulations is
\begin{equation}\label{3marked}
|G^3(n,n-3)|=(n-2)\G{n}{n-3}.
\end{equation}

\begin{definition}
Let $\Phi\in G(n)$ and consider the representation of $\Phi$ as a set of points in the interior or boundary of a regular $n$-gon embedded in $\mathbb R^2$ and centered at the origin. Any component of $\Phi$ is  a subset of the planar region. There is thus a unique component of $\Phi$ containing the origin. We call this component the
 {\em central polygon} $Z(\Phi)$ of a dissection $\Phi$. Let $G_m(X)$ be the subset of $G(X)$ consisting of dissections whose central polygon is an $m$-gon; \[ G_m(X)=\{\Phi\in G(X): Z(\Phi)\in \PP_m(\Phi)\},\] and put $G_{\ne m}(X)=G(X)\setminus G_m(X)$.
\end{definition}

For $\Phi\in G_m(n,k)$, if the regions outside of $Z(\Phi)$ are triangulated then $m=n-k$. More generally,
\begin{equation}\label{mn-k}
m\le n-k \text{ for any $\Phi\in G_m(n,k)$}.
\end{equation}

\begin{definition}
Let $\Phi\in G(n)$. Given an edge $xy\in E(\Phi)$ and a vertex $v\in V(\Phi)$, we say that $v$ is {\em outer to} $xy$ if $v$ lies strictly between $x$ and $y$ on the shorter path of the $n$-gon connecting them.
\end{definition}

\begin{remark}\label{outer}
A vertex $v\in V(Z(\Phi))$ cannot be outer to any edge $xy\in E(\Phi)$.
\end{remark}

\section{Rotationally symmetric dissections}\label{rotsec}

The enumeration of the sets $G(n,k; \R^{n/d})$ of rotationally symmetric dissections can be achieved by considering separately the following two classes.
 
 \begin{definition}\label{classes}
 A dissection $\Phi\in G(n,k; \R^{n/d})$ is said to be {\em centrally bordered} if $\Phi\in G_d(n,k; \R^{n/d})$ and {\em centrally undbordered} if
 $\Phi\in G_{\ne d}(n,k; \R^{n/d})$. 
 \end{definition}
 Lemma \ref{bdd} addresses the case of centrally bordered dissections; Lemmas \ref{Frange}, \ref{FU} and \ref{UF} and  Theorem \ref{bijthm} address the case of centrally unbordered dissections. When considering the set  $G(n; \R^{n/d})$ it is convenient to put $j=n/d$.
Let $\delta_{xy}$ denote the Kronecker delta.

\begin{lemma} \label{bdd}
Let $d,j\ge 2$ and let $0 \le k\le n-3$.
\begin{equation}\label{bddeq}
|G_d(n,k;\R^j)|= j \G{j+1}{(k-d+\delta_{d2})/d} \ .
\end{equation}
\begin{proof}
    Let $\Phi \in G_d(n,k;\R^j)$. By symmetry the central polygon $Z(\Phi)$ is a regular $d$-gon, which can be positioned in $j$ different ways in the $n$-gon. Since all $d-\delta_{d2}$ edges of $Z(\Phi)$ are diagonals of $\Phi$, there are $k-d+\delta_{d2}$ diagonals of $\Phi$ which are not edges of $Z(\Phi)$. By symmetry these diagonals are equally distributed among the $d$ resulting $(j+1)$-gons, giving
$\G{j+1}{(k-d+\delta_{d2})/d}$ choices for each position of $Z(\Phi)$.
\end{proof}
\end{lemma}

\begin{figure}
\centering

\subfigure[$G_{4}(12,6; \R^6)\leftrightarrow G^2(6,2)$]
{
   \includegraphics[scale =0.8] {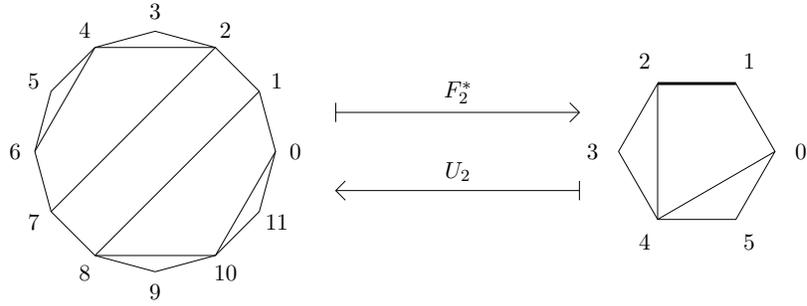}
   \label{fig:subfig1}
 }

 \subfigure[$G_{4}(12,6; \R^6)\leftrightarrow G^2(6,2)$]
 {
   \includegraphics[scale =0.8] {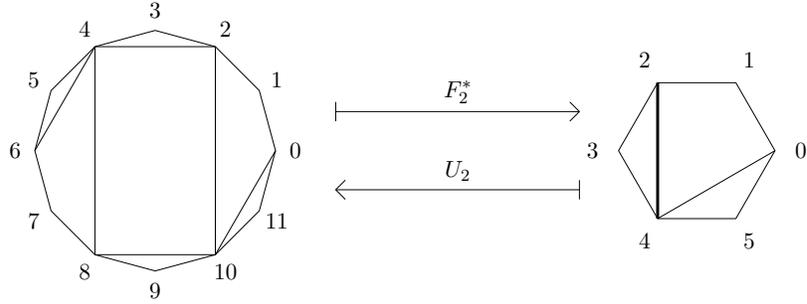}
   \label{fig:subfig2}
 }

 \subfigure[$G_{6}(12,6; \R^6)\leftrightarrow G^3(6,3)$]
 {
   \includegraphics[scale =0.8] {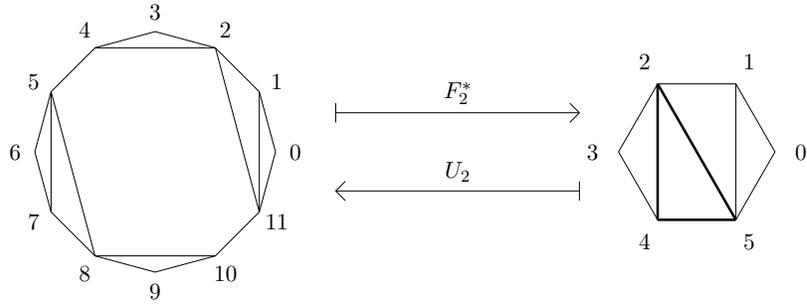}
   \label{fig:subfig3}
 }

 \subfigure[$G_{6}(12,6; \R^4)\leftrightarrow G^2(4,1)$]
 {
   \includegraphics[scale =0.8] {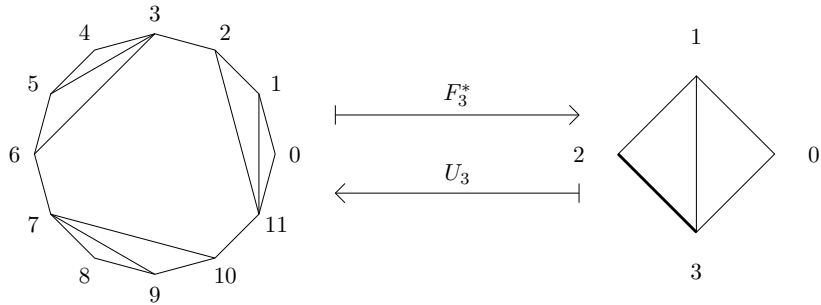}
   \label{fig:subfig4}
 }

\caption{Examples showing the bijections of Theorem \ref{bijthm}. \label{bijfig}}
\end{figure}

Enumeration of the centrally unbordered dissections $G_{\ne d}(n,k;\R^{n/d})$ is accomplished by the introduction of a bijection with the marked dissections which were enumerated in Lemma \ref{marked}.
We define the following ``furling" maps $F_d$ and $F^*_d$ (see Figure \ref{bijfig}).

\begin{definition}
Let $d,j\ge 2$.
\begin{enumerate}
\item Define a function $f_d$ on edges of a graph of $n$ vertices by $f_d(xy)=[x]_j[y]_j$.
\item Let $\Phi$ be a dissection or a component of a dissection of an $n$-gon. Define the labeled graph $F_d(\Phi)$ by
$V(F_d(\Phi))=\{[x]_j : x\in V(\Phi)\}$ and $E(F_d(\Phi))=f_d[E(\Phi)]$. \footnote{Functions on subsets of a set are defined in the usual way and denoted using square brackets.}
\item For $\Phi\in G_{\ne d}(n,k;\R^j)$, define  $F^*_d(\Phi)=(F_d(\Phi), F_d(Z(\Phi))).$\end{enumerate}
\end{definition}

The conclusions of the following remark are easy observations.
\begin{remark}
Let $d,j\ge 2$.
\begin{enumerate}[(i)]
\item If $\Phi\in G_{\ne d}(n,k;\R^j)$ then its central polygon $Z(\Phi)$ is itself invariant under $\R^j$ and hence $Z(\Phi)$ is an $rd$-gon for some $r\ge 2$. Therefore
$G_{\ne d}(n,k;\R^j)$ can be partitioned as follows.
\begin{equation}\label{partition}
G_{\ne d}(n,k;\R^j)=\bigcup_{r\ge 2}G_{rd}(n,k;\R^j).
\end{equation}
\item  If $\Phi \in G_{\ne d}(n,k; \R^j)$ and $xy\in E(\Phi)$, then  the distance between $x$ and $y$ is at most $j-1$.
\item Suppose $\Phi \in G_{\ne d}(n,k; \R^j)$. From the previous observation and by symmetry, it follows that for $0\le x<y\le j-1$,
\begin{equation}\label{xy}
xy\in E(F_d(\Phi))  \text{ if and only if either } xy \in E(\Phi) \text{ or } y(x+j)\in E(\Phi).
\end{equation}
\end{enumerate}
\end{remark}


The map ${\bf u}_d$ will be used to output symmetrically distributed edges in a dissection.

\begin{definition}
Let $d,j\ge 2$. For an edge $ab$ of a dissection  $\Phi\in G(j)$, define \[{\bf u}_d(ab)=\bigcup\limits_{0\le i \le d-1} \Big \{ \{[a+ij]_n,[b+ij]_n\}\Big \} .\]
\end{definition}

 Note that if  $\Phi\in G_{rd}(n,k; \R^{j})$ and $\gamma=Z(\Phi)$ then  by the $d$-fold symmetry, for some $0\le v_r=v_0<\ldots <v_{r-1}\le j-1$,
\begin{align}\label{symmcenter}
V(\gamma)=\{v_0,\ldots, v_{r-1},v_0+j,\ldots, v_{r-1}+j, \ldots, v_{r-1}+(d-1)j\}, \nonumber\\
\nonumber \\
E(\gamma) = {\bf u}_d(v_0v_1) \cup {\bf u}_d(v_1v_2) \cup \ldots \cup {\bf u}_d(v_{r-2}v_{r-1})\cup {\bf u}_d(v_{r-1}(v_0+j))
\nonumber \\
\end{align}

The following lemma is readily verified using \eqref{xy} and Remark \ref{outer}. It employs the notation \eqref{symmcenter} for $Z(\Phi)$.
\begin{lemma}\label{preimage}
Let  $d,j\ge 2$ and let $\Phi\in G_{rd}(n,k; \R^j)$. 
Consider the map $f_d:E(\Phi)\to E(F_d(\Phi))$ defined above. For an edge $ab\in E(F_d(\Phi))$, with $a<b$,  the preimage of $ab$ under $f_d$ is given by:
\begin{align}\label{preimage eq}
f_d^{-1}[ab]=
\begin{cases}
{\bf u}_d(ab) , \qquad \text{ if $a>v_0$ or $b<v_{r-1}$},\\
{\bf u}_d(b(a+j)), \qquad \text{ if $a\le v_0$,  $b\ge v_{r-1}$, and $ab\ne v_0v_1$}, \\
{\bf u}_d(ab)\cup {\bf u}_d(b(a+j)), \qquad \text { if $ab=v_0v_1$.}\\
\end{cases}
\end{align}
\end{lemma}

\begin{lemma}\label{Frange}
Let $d, j, r\ge 2$ and let $1\le k\le n-3$. If $\Phi \in G_{rd}(n,k; \R^j)$ then
 \[F^*_d(\Phi)\in {G^r(j,k/d-\delta_{r2}}).\]
\begin{proof}
Clearly $F_d(\Phi)$ has $j$ vertices. We show that its diagonals are noncrossing. Suppose that $a_1b_1$ and $a_2b_2$ are crossing diagonals of $F_d(\Phi)$, with $0\le a_1<a_2<b_1<b_2\le j-1$. By Lemma \ref{preimage}, for each $i=1, 2$ either $a_ib_i$ or $a_i(b_i+j)$ is a diagonal of $\Phi$. Since $a_1<a_2<b_1<b_2<a_1+j<a_2+j$, these two diagonals of $\Phi$ are crossing. This contradiction shows that $F_d(\Phi)\in G(j)$.

We next show that $F_d(Z(\Phi))$ is an $r$-component of $F_d(\Phi)$. By symmetry the center $Z(\Phi)$ has the form \eqref{symmcenter}. Therefore $V(F_d(Z(\Phi))=\{v_0,\ldots, v_{r-1}\}$ and
$E(F_d(Z(\Phi)))=\{v_0v_1, \ldots, v_{r-1}v_r\}$.
Now if $v_sv_t \in E(F_d(\Phi))$ with $v_s<v_t$ then by \eqref{xy} either $v_sv_t$ or $v_t(v_s+j)$ is in $E(\Phi)$. Therefore by the fact that $Z(\Phi)$ is a component of $\Phi$ and by \eqref{component}, either $v_sv_t$ or $v_t(v_s+j)$ is in $E(Z(\Phi))$. 
Thus $v_sv_t\in E(F_d(Z(\Phi)))$. Applying \eqref{component} again gives the conclusion.

Let $l$ be the number of diagonals of $F_d(\Phi)$.
Suppose $r\ge 3$. In this case an edge $xy \in E(\Phi)$ is a diagonal of $\Phi$ if and only if $f_d(xy)$ is a diagonal of $F_d(\Phi)$. By \eqref{preimage eq}, for each diagonal
$ab$ of $F_d(\Phi)$, the preimage $f_d^{-1}[ab]$ consists of $d$ diagonals of $\Phi$. Furthermore, if $ab\ne a'b'$ then $f_d^{-1}[ab]$ and $f_d^{-1}[a'b']$ are disjoint. Thus $k=dl$.
Now suppose $r=2$. If $ab$ is a diagonal of $F_d(\Phi)$ with $ab\ne v_0v_1$ then $f_d^{-1}[ab]$ consists of $d$ diagonals of $\Phi$.
Note that either $v_0v_1$ or $v_1(v_0+j)$ is a diagonal of $\Phi$, since otherwise $k=0$. If both $v_0v_1$ and $v_1(v_0+j)$ are diagonals then
$f_d^{-1}[v_0v_1]$ consists of $2d$ diagonals of $\Phi$. If only one of $v_0v_1$ and $v_1(v_0+j)$ is a diagonal of $\Phi$ then $v_0v_1$ is not a diagonal of $F_d(\Phi)$ and $f_d^{-1}[v_0v_1]$ consists of $d$ diagonals of $\Phi$.  Thus in either case for $r=2$, $k=dl+d$ and the result follows.
\end{proof}
\end{lemma}

Lemma \ref{Frange} shows that  $F^*_d:G_{\ne d}(n,k; \R^j)\to G^*(j)$. We next define an ``unfurling" map $U_d$.
\begin{definition}\label{Ud def}
Let $d,j\ge 2$. Define $U_d: G^*(j)\rightarrow G(n)$ as follows.
Let $(\Theta, \beta)\in G^{r}(j)$,  and denote the vertices of $\beta$ by $v_0<\ldots <v_{r-1}$.
Define $U_d(\Theta,\beta)$ by $V(U_d(\Theta,\beta))=\{0,\ldots,n-1\}$ and $E(U_d(\Theta,\beta))=\bigcup\limits_{ab\in E(\Phi)} f_d^{-1}[ab]$, where $f_d^{-1}$ is given by \eqref{preimage eq}.
\end{definition}

\begin{lemma}\label{FU}
Let $d, j \ge 2$. If $(\Theta,\beta)\in G^*(j)$ then $F^*_d(U_d(\Theta,\beta))=(\Theta,\beta)$.
\begin{proof}
Clearly $V(F_d(U_d(\Theta,\beta)))=V(\Theta)$ and
\begin{align*}
E(F_d(U_d(\Theta, \beta)))=f_d[E(U_d(\Theta,\beta))]=f_d \left [\bigcup\limits_{ab\in E(\Theta)} f_d^{-1}[ab] \right ] =E(\Theta),
\end{align*}
so $F_d(U_d(\Theta, \beta))=\Theta$.
Suppose $(\Theta, \beta)\in G^r(j)$; denote the vertices of $\beta$ by $v_r=v_0<  \ldots <v_{r-1}$, and let $\gamma$ be the graph given by \eqref{symmcenter}.
By definition $\gamma$ is a subgraph of $U_d(\Theta, \beta)$.
As in the the proof of Lemma \ref{Frange}, the condition \eqref{component} can be used to show that in fact $\gamma$ is a component of $U_d(\Theta, \beta)$.
Finally since the vertices of $\gamma$ include the regular $d$-gon with vertices $v_0, v_0+j, \ldots , v_0+(d-1)j$, their convex hull contains the origin, so
$\gamma=Z(U_d(\Theta,\beta))$.
Thus $F_d(Z(U_d(\Theta, \beta)))= F_d(\gamma)=\beta$, completing the proof.
\end{proof}
\end{lemma}

\begin{lemma}\label{UF}

Let $d,j\ge 2$. If $\Phi\in G_{\ne d}(n;\R^j)$ then $U_d(F_d^*(\Phi))=\Phi$.
\begin{proof}
Let $\Phi\in G_{rd}(n;\R^j)$. It is easily seen that $V(U_d(F_d^*(\Phi)))=V(\Phi)$.
As above, the center $Z(\Phi)$ has the form \eqref{symmcenter}.
 Therefore $V(F_d(Z(\Phi)))=\{v_0,\ldots, v_{r-1}\}$, and 
by Definition \ref{Ud def},
\begin{align*}
E(U_d(F_d^*(\Phi)))=E(U_d(F_d(\Phi), F_d(Z(\Phi)))
				=\bigcup\limits_{ ab\in E(F_d(\Phi))} f_d^{-1}[ab]
				=E(\Phi).
\end{align*}

\end{proof}
\end{lemma}

\begin{theorem}\label{bijthm}
Let $1\le k \le n-3$, let $r\ge 2$ and let $2\le d\le n/3$ with $d | n$. Then there exists a bijection:
\begin{equation*}\label{bijeq}
G_{rd}(n,k;\R^j) \longleftrightarrow G^r(j,k/d-\delta_{r2}).
\end{equation*}
\begin{proof}
By \eqref{partition} and Lemmas \ref{Frange}, \ref{FU} and \ref{UF}, the bijection is given in one direction by $F_d^*$  and in the other direction by $U_d$.
\end{proof}
\end{theorem}

Lemma \ref{marked} and Theorem \ref{bijthm} imply that for $d\le n/3$,
\begin{align}\label{ubeq}
|G_{\ne d}(n,k;\R^{n/d})|
&=|G^2(n/d,k/d-1)|+\summ_{r\ge 3}{|G^r(n/d,k/d)|}\nonumber \\
&=
{n+k-d\over d}\G{n/d}{k/d-1}+\summ_{\substack{r\ge 3; \ n_1+\ldots+n_r=n/d \nonumber \\
k_1+\ldots+k_r +|\{i|n_i\ge 2\}|= k/d}} {n \over r}  \prodd_{i=1}^r\G{n_i+1}{k_i}.\\
\end{align}
Note that if $d> n/3$ then $|G_{\ne d}(n,k;\R^{n/d})|=0$.

\subsection{Proof of Theorems \ref{dih} and \ref{cyc}}

The proof of Theorems  \ref{dih} and \ref{cyc} now follows by substituting into equations \eqref{dihgeneq2} and \eqref{cycgeneq2} the expressions obtained in
\eqref{p=0}
--\eqref{Todd} for the number of axially symmetric dissections, and the values obtained in \eqref{bddeq} and \eqref{ubeq} for rotationally symmetric dissections.

\section{Interesting special cases}\label{specs}

The enumeration formulas  can be specialized to certain classes of dissections, namely for specific values of $n-k$ and for specific values of $k$.
The next lemma is equivalent to Catalan's $k$-fold convolution formula \cite{Cat,LF}.

\begin{lemma}\label{catcon}
For any $n,m\ge 0$,
\begin{equation*}\label{catid}
\summ_{\substack{i_1+\ldots+i_m=n\\i_1,\ldots,i_m\ge 0}}C_{i_1}\cdots C_{i_m}=
\begin{cases}
\frac{m(n+1)(n+2)\cdots(n+\frac m2-1)}{2(n+\frac{m}{2}+2)(n+\frac{m}{2}+3)\cdots (n+m)} C_{n+m/2}, \qquad &\text{if $m$ is even}\\
\\
\frac{m(n+1)(n+2)\cdots(n+\frac{m-1}2)}{(n+\frac{m+3}{2})(n+\frac{m+3}2+1)\cdots(n+m)} C_{n+(m-1)/2}, \qquad &\text{if $m$ is odd.}
\end{cases}
\end{equation*}
\end{lemma}

\begin{lemma}\label{disscons}
\begin{enumerate}
\item For any $n\ge 2$,
\begin{equation}\label{NCT}
\G{n}{n-3}+\G{n}{n-2}=C_{n-2}.
\end{equation}
\item For any $n\ge 2, q\ge 2$,
\begin{equation} \label{delta}
\summ_{i+j=n}\G{i+1}{i-1} \G{j+1}{j+1-q}= \G{n}{n-q}.
\end{equation}
\item For any $n\ge 3$,
\begin{equation}\label{disscon1}
\summ_{i+j=n}\G{i+1}{i-2} \G{j+1}{j-2}=C_{n-1}-2C_{n-2}.
\end{equation}
\item For any $n\ge 3$,
\begin{equation}\label{disscon2}  
\summ_{i+j=n}\G{i+1}{i-2} \G{j+1}{j-3}= {(n-3)(n-4)\over 2n}C_{n-2}.
\end{equation}
\end{enumerate}
\begin{proof}
Equations \eqref{NCT} and \eqref{delta} follow from \eqref{catdiss}, and \eqref{disscon1} follows from \eqref{NCT} and from Lemma \ref{catcon}.
To prove \eqref{disscon2}, we show that
\begin{equation}\label{disconn2a}
(n-4)\G{n}{n-4}=n\summ_{i+j=n}\G{i+1}{i-2} \G{j+1}{j-3}.
\end{equation}
The result will then follow since $\G{n}{n-4}={n-3\over 2}C_{n-2}$. Now the left hand side of \eqref{disconn2a} is the number
of almost-triangulations marked by a diagonal (i.e., $(\Phi, \beta)\in G^2(n,n-4)$ where $V(\beta)$ is not of the form $\{v,[v+1]_n\}$).
These can also be enumerated as follows. Choose one vertex $v$ out of the $n$ vertices,
then choose $2\le i\le n-3$ and $j=n-i$. Mark the diagonal $v[v+i]_n$, and choose a triangulation of the resulting $(i+1)$-gon and an almost-triangulation of the resulting $(j+1)$-gon.
\end{proof}
\end{lemma}

The details of the proof of Theorem \ref{k=n-5} for $n$ even are given below. Most of the details of the other cases are omitted.
Note that if
\[\summ_{\substack{\ n_1+\ldots+n_r=n/d\\ k_1+\ldots+k_r + |\{i : n_i\ge 2\}| = k/d}}\G{n_i+1}{k_i}\ne 0\] for some $d\ge 2$, $r\ge 3$, then it follows from \eqref{mn-k}
that $6\le rd\le n-k$. Thus these terms vanish in the cases $k=n-3$, $n-4$ or $n-5$.

Applying \eqref{diheq} and \eqref{cyceq} when $k=n-3$ (i.e., for triangulations), recovers the result of Moon and Moser %
 and the result of Brown. (We omit the details here as the even dihedral case of Theorem \ref{k=n-5}, for which the details are provided, is a similar calculation).
\begin{theorem}\label{k=n-3}
\begin{enumerate}
\item {\cite{MM}} Let $n\ge 3$. The number of triangulations of an $n$-gon (equivalently, the number of vertices of the $(n-3)$-dimensional associahedron) modulo the dihedral action is
\begin{equation*}
|G(n,n-3)/D_{2n}| =
\begin{cases}
\frac{1}{2n}C_{n-2} + \frac{1}{3}C_{n/3-1}+\frac{3}{4}C_{n/2-1}, \qquad \qquad &\text{if $n$ is even} \\
\frac{1}{2n}C_{n-2} + \frac{1}{3}C_{n/3-1}+\frac{1}{2}C_{(n-3)/2}, \qquad \qquad &\text{if $n$ is odd.}
\end{cases}
\end{equation*}
\item {\cite{B}} Let $n\ge 3$. The number of triangulations of an $n$-gon (equivalently, the number of vertices of the $(n-3)$-dimensional associahedron) modulo the cyclic action is \[|G(n,n-3)/Z_n|=\frac{1}{n}C_{n-2}+\frac{1}{2}C_{n/2-1}+\frac{2}{3}C_{n/3-1}.\]
\end{enumerate}
\end{theorem}


%
Setting $k=n-4$ in \eqref{diheq} recovers the following result of the authors.
\begin{theorem}\label{ges}\cite{BR}
Let $n\ge 4$, and let $g^{(e)}(n)$ be the number of almost-triangulations of an $n$-gon (equivalently, edges of the $(n-3)$-dimensional associahedron) modulo the dihedral action. Then
\begin{equation*}
g^{(e)}(n)=
\begin{cases}
 (\frac{1}{4}-\frac{3}{4n})C_{n-2}+\frac{3}{8}C_{n/2-1}+(1-\frac{3}n)C_{n/2-2}+\frac{1}{4}C_{n/4-1}, &  \text{ if $n$ is even}\\
 (\frac{1}{4}-\frac{3}{4n})C_{n-2}+\frac{1}{4}C_{(n-3)/2}, & \text{ if $n$ is odd}. \\
\end{cases}
\end{equation*}
\end{theorem}
Setting $k=n-4$ in \eqref{cyceq} gives the result of Theorem \ref{cycn4}.
The following proposition gives the details needed to simplify \eqref{diheq}, thus completing the proof of the even dihedral case of Theorem  \ref{k=n-5}. The other cases are easier.

\begin{proposition}\label{n-5prop}
Let $n\ge 6$ and $k=n-5$. If $n$ is even then
\begin{equation}\label{Gn5E}
{1\over 2n} \G{n}{k}={(n-3)^2(n-4)\over 8n(2n-5)}C_{n-2},
\end{equation}
\begin{equation}\label{Gn5T0D}
{1\over 2} \G{n/2+1}{(k-1)/2}={n-4\over 8}C_{n/2-1},
\end{equation}
\begin{equation}\label{Gn5T0}
{1\over 4} \G{n/2+1}{k/2} =0,
\end{equation}
\begin{equation}\label{Gn5bdd}
\summ_{3\le d|n}{\varphi(d)\over 2d} \G{n/d+1}{k/d-1} =  {2\over 5}C_{n/5-1},
\end{equation}
\begin{equation}\label{Gn5ub2}
\summ_{2\le d\le n/3}{\varphi(d)(n+k-d) \over 2dn}\G{n/d}{k/d-1}  =0,
\end{equation}
\begin{equation}\label{Gn5ub}
\summ_{\substack{2\le d|n; \ r\ge 3; \ n_1+\ldots+n_r=n/d\\ k_1+ \ldots+k_r + |\{i : n_i\ge 2\}|= k/d}} {\varphi(d)\over 2r}  \prodd_{i=1}^r\G{n_i+1}{k_i} =0,
\end{equation}
\begin{equation}\label{Gn5T}
{1\over 4} \summ_{\substack{1\le t \le k\\ n_0+\ldots + n_t=n/2\\k_0+\ldots + k_t=(k-t)/2}}
\prodd_{s=0}^ {t} \big (\G{n_s+1}{k_s-1} + \G{n_s+1}{k_s} \big)  = {n^2-2n-12\over 16(n-3)}C_{n/2-1},
\end{equation}
and
\begin{equation}\label{Gn5TR}
{1\over 4} \summ_{\substack{0\le t \le k\\ n_0+\ldots + n_t=n/2-1\\ k_0+\ldots + k_t=(k-t)/2}} \prodd_{s=0}^ {t} \big (\G{n_s+1}{k_s-1} + \G{n_s+1}{k_s} \big)  = {n\over 16(n-3)} C_{n/2-1}.
\end{equation}
\begin{proof}
Equation \eqref{Gn5E} follows from \eqref{catdiss}, and \eqref{Gn5T0D} and \eqref{Gn5T0} are immediate.
In \eqref{Gn5bdd}, the only nonzero summand corresponds to $d=5$.
To prove \eqref{Gn5ub2}, note that if $d$ divides $n$ and $k$ then $d=5$, but then $\G{n/d}{k/d-1}=0$ since $n/d\ge 3$.
Equation \eqref{Gn5ub} follows from the remark before Theorem \ref{k=n-3}.
For \eqref{Gn5T}, note that if $\prodd_{s=0}^p (\G{n_s+1}{k_s-1} + \G{n_s+1}{k_s})\ne 0$, then $k_s\le n_s-1$ for all $s$. Therefore in this case
\[(k-t)/2=\summ_{s=0}^t k_s \le \summ_{s=0}^t(n_s-1)=n/2-t, \]
i.e., $t\le n-k-2$, with equality if and only if $k_s=n_s-1$ for all $s$.
Therefore, using the notation of Section \ref{axsec}, the nonzero summands in \eqref{Gn5T} correspond to $|G(n,n-5;\T;1)|$ and $|G(n,n-5;\T;3)|$. Now
\begin{align}\label{eqn}
|G(n,n-5;\T;1)|&=\summ_{\substack{n_0+ n_1=n/2\\ k_0 + k_1=n/2-3}} \left (\G{n_0+1}{k_0-1} + \G{n_0+1}{k_0}\right)\left(\G{n_1+1}{k_1-1} + \G{n_1+1}{k_1}\right ) \qquad \qquad \nonumber \\
    &=\summ_{\substack{n_0+ n_1=n/2\\ k_0 + k_1=n/2-3}}
    \big (\G{n_0+1}{k_0-1} \G{n_1+1}{k_1-1}+\G{n_0+1}{k_0-1} \G{n_1+1}{k_1}
    + \G{n_0+1}{k_0} \G{n_1+1}{k_1-1} + \G{n_0+1}{k_0} \G{n_1+1}{k_1}\big ).\nonumber\\
\end{align}
Any nonzero terms in \eqref{eqn} have $(k_0,k_1)=(n_0-1,n_1-2)$ or $(k_0,k_1)=(n_0-2,n_1-1)$. Therefore by Lemma \ref{disscons} and by symmetry,
\begin{align*}
|G(n,n-5;\T;1)|=
    2\summ_{n_0+ n_1=n/2}
    \big (\G{n_0+1}{n_0-2} \G{n_1+1}{n_1-3}+ \G{n_0+1}{n_0-2} \G{n_1+1}{n_1-2} + \G{n_0+1}{n_0-1} \G{n_1+1}{n_1-3} + \G{n_0+1}{n_0-1} \G{n_1+1}{n_1-2}\big )\\
    =2\left [{(n-6)(n-8)\over 4n} C_{n/2-2} + C_{n/2-1}-2C_{n/2-2} + {n-6\over 4}C_{n/2-2} + C_{n/2-2}\right ].
\end{align*}
Next by Lemma \ref{catcon},
\begin{multline*}
|G(n,n-5;\T;3)|=\\
\summ_{n_0+n_1+n_2+ n_3=n/2}\left (\G{n_0+1}{n_0-2}+\G{n_0+1}{n_0-1}\right ) \left (\G{n_1+1}{n_1-2}+\G{n_0+1}{n_0-1}\right )
 \left (\G{n_2+1}{n_2-2}+\G{n_0+1}{n_0-1}\right )\left (\G{n_3+1}{n_3-2}+\G{n_0+1}{n_0-1}\right )\\
=\summ_{n_0+n_1+n_2+ n_3=n/2}C_{n_0-1}C_{n_1-1}C_{n_2-1}C_{n_3-1}={2n-12\over  n}C_{n/2-2}.
\end{multline*}
Equation \eqref{Gn5T} now follows by simplifying these expressions and using the relation $C_{n/2-2}={n\over 4(n-3)}C_{n/2-1}$. A similar argument proves \eqref{Gn5TR}.
\end{proof}

\end{proposition}

The proof of Theorem \ref{k=n-5} now follows by applying Proposition \ref{n-5prop} and Theorems \ref{dih} and \ref{cyc}.
The proof of Theorems \ref{k=n-6} and \ref{k=2} proceeds along similar lines. For Theorem \ref{k=n-6}, note that if $|G^r(n/d,k/d-\delta_{r2})|\ne 0$ then either $r=2$ and $d=2$; or $r=2$ and $d=3$; or $r=3$ and $d=2$ (see Figure \ref{bijfig}). The last of these cases can be computed using \eqref{3marked}:
\[|G^3(n/2,(n-6)/2-3+3)|
=(n/2-2)\G{n/2}{n/2-3}=(n/2-2)C_{n/2-2}.\]

\end{document}